\documentclass[a4paper,12pt]{article}
\usepackage{amssymb}
\usepackage{amsthm}
\usepackage{amsxtra}
\usepackage{amsfonts}
\usepackage{xr}

\begin{document}

\def\sect{\section}

\newtheorem{thm}{Theorem}[section]
\newtheorem{cor}[thm]{Corollary}
\newtheorem{lem}[thm]{Lemma}
\newtheorem{prop}[thm]{Proposition}
\newtheorem{propconstr}[thm]{Proposition-Construction}

\theoremstyle{definition}
\newtheorem{para}[thm]{}
\newtheorem{ax}[thm]{Axiom}
\newtheorem{conj}[thm]{Conjecture}
\newtheorem{defn}[thm]{Definition}
\newtheorem{defns}[thm]{Definitions}
\newtheorem{notation}[thm]{Notation}
\newtheorem{rem}[thm]{Remarks}
\newtheorem{remark}[thm]{Remark}
\newtheorem{question}[thm]{Question}
\newtheorem{example}[thm]{Example}
\newtheorem{problem}[thm]{Problem}
\newtheorem{excercise}[thm]{Exercise}
\newtheorem{ex}[thm]{Exercise}

\def\Bbb{\mathbb}
\def\cal{\mathcal}
\def\mL{{\mathcal L}}
\def\mC{{\mathcal C}}

\overfullrule=0pt

\def\si{\sigma}
\def\prf{\smallskip\noindent{\it        Proof}. }
\def\call{{\cal L}}
\def\nat{{\Bbb  N}}
\def\la{\langle}
\def\ra{\rangle}
\def\inv{^{-1}}
\def\ld{{\rm    ld}}
\def\trdeg{{tr.deg}}
\def\dim{{\rm   dim}}
\def\th{{\rm    Th}}
\def\rest{{\lower       .25     em      \hbox{$\vert$}}}
\def\ch{{\rm    char}}
\def\zee{{\Bbb  Z}}
\def\conc{^\frown}
\def\acl{{\rm acl}}
\def\cls{cl_\si}
\def\cals{{\cal S}}
\def\mult{{\rm  Mult}}
\def\calv{{\cal V}}
\def\aut{{\rm   Aut}}
\def\ffi{{\Bbb  F}}
\def\ffiti{\tilde{\Bbb          F}}
\def\degs{deg_\si}
\def\calx{{\cal X}}
\def\gal{{\rm Gal}}
\def\cl{{\rm cl}}
\def\loc{{\rm locus}}
\def\calg{{\cal G}}
\def\calq{{\cal Q}}
\def\calr{{\cal R}}
\def\caly{{\cal Y}}
\def\aff{{\Bbb A}}
\def\cali{{\cal I}}
\def\calu{{\cal U}}
\def\epsilon{\varepsilon} 
\def\Uu{{\cal U}}
\def\rat{{\Bbb Q}}
\def\ga{{\Bbb G}_a}
\def\gm{{\Bbb G}_m}
\def\cee{{\Bbb C}}
\def\ree{{\Bbb R}}
\def\frob{{\rm Frob}}
\def\Frob{{\rm Frob}}
\def\fix{{\rm Fix}}
\def\Uu{{\cal U}}
\def\proj{{\Bbb P}}
\def\sym{{\rm Sym}}
\def\Ker{{\rm Ker}}
\def\dcl{{\rm dcl}}
\def\calm{{\mathcal M}}

\font\helpp=cmsy5
\def\semdp
{\hbox{$\times\kern-.23em\lower-.1em\hbox{\helpp\char'152}$}\,}

\def\dnfo{\,\raise.2em\hbox{$\,\mathrel|\kern-.9em\lower.35em\hbox{$\smile$}
$}}
\def\dnf#1{\lower1em\hbox{$\buildrel\dnfo\over{\scriptstyle #1}$}}
\def\dfo{\;\raise.2em\hbox{$\mathrel|\kern-.9em\lower.35em\hbox{$\smile$}
\kern-.7em\hbox{\char'57}$}\;}
\def\df#1{\lower1em\hbox{$\buildrel\dfo\over{\scriptstyle #1}$}}        
\def\stab{{\rm Stab}}
\def\qfcb{\hbox{qf-Cb}}
\def\perf{^{\rm perf}}
\def\sipm{\si^{\pm 1}}

\def\vlabel{\label}

\title{Geometric representation in the theory of pseudo-finite fields}

\author{\"Ozlem Beyarslan\thanks{partially supported by Bo\u{g}azi\c{c}i BAP
  14B06P1-7840}\ \footnotemark[3]
\and 
Zo\'e Chatzidakis\thanks{partially supported by PITN-2009-238381
    and by ANR-13-BS01-0006}\ \thanks{Both authors would like to thank the Mathematical Sciences Research
Institute in Berkeley for its hospitality and support during the spring
2014, when part of this research was conducted. 
Research at MSRI was supported by
NSF grant DMS-0932078 000}
}
\date{}
%\centerline{\today}
\maketitle

\begin{abstract} We study the automorphism group of the algebraic
  closure of a substructure $A$ of a pseudo-finite field $F$, or more
  generally, of a bounded  PAC field $F$. This paper answers some
  of the questions of \cite{BH}, and in particular that any finite group
  which is geometrically represented in a pseudo-finite field must be abelian. 
\end{abstract}

\section*{Introduction}

\noindent This paper investigates the relationship between model-theoretic
definable closure and model-theoretic algebraic closure in certain
fields. In other words: if $F$ is a field, and $A\subseteq F$ satisfies
$A=\dcl(A)$, what can one say of the group $\aut(\acl(A)/A)$ of
restrictions to $\acl(A)$ of elements of $\aut(F/A)$? When is it
non-trivial? A natural assumption to add is to look at a slightly
smaller group, and to impose on $A$ that it contains an elementary
substructure of $F$. Indeed, we certainly want to impose that our
automorphisms fix $\acl^{eq}(\emptyset)$. \\
This paper extends some of the results of \cite{BH}, with
completely new
proofs, and answers some of the questions there. Here are the main results we
obtain:\\[0.1in]
{\bf Theorem \ref{thm1}}. {\em Let $F$ be a bounded field, $A=\dcl(A)$ a
  subfield of $F$ containing an elementary substructure of $F$, and let $p$ be a prime dividing
  $\#(\aut(\acl(A)/A)$ and $\#G(F)$. Then $p\neq {\rm char}(F)$,
  and $\mu_{p^\infty}\subset F(\zeta_p)$.}\\[0.1in]
{\bf Theorem \ref{thm2}}. {\em Let $F$ be a pseudo-finite field, [or more
  generally a bounded PAC field].  
Assume that for some  subfield  $A=\dcl(A)$ 
 of $F$ containing an elementary substructure of $F$, the group
 $G:=\aut(\acl(A)/A)$ is non-trivial. [Assume in addition
that all primes dividing $\#G$ divide $\#G(F)$]. \\
Then $G$ is abelian, for any prime $p$ dividing $\#G$, we have $p\neq
{\rm char}(F)$, and 
$\mu_{p^\infty}\subset F$.}\\[0.1in]
We  give an example (\ref{prop1}) which shows that the hypotheses on
$F_0$ cannot be weakened to assume that $A$ contains a substucture $F_0$
with $\acl^{eq}(\emptyset)\subset \dcl^{eq}(F_0)$. We also give a
partial answer to a question of \cite{BH} on centralisers. \\[0.1in]

\section{The results}

\begin{notation} Let $F$ be a field. Throughout the paper, $\dcl$ and
  $\acl$ will denote the model-theoretic  definable  and algebraic
  closures, taken within the structure $F$ or possibly some elementary extension
  of $F$. \\
$F^{alg}$ denotes an algebraic closure of $F$ (i.e., an algebraically
closed field containing $F$ and minimal such), $F^s$ its separable
closure  and $G(F)$ its absolute
Galois group, i.e., $\gal(F^{s}/F)$. \\
If $A\subset B$ are subfields of $F$, we denote by $\aut(B/A)$ the set
of automorphisms of $B$ which preserve all $\call(A)$-formulas true in
$F$, and by $\aut_{\rm field}(B/A)$ the set of (field) automorphisms of
$B$ which fix the elements of $A$. \\ 
We let ${\mathbb\mu}_{p^\infty}$ denote the group of all
  $p^n$-th roots of unity if $p\neq {\rm char}(F)$, and $\zeta_p$ a
  primitive $p$-th root of unity.\\
Let $G_1, G_2$ be  profinite groups, $p$ a prime. We say that $p$ divides
  $\#G_1$ if $G_1$ has a finite quotient with order divisible by $p$. We
  write $(\#G_1,\#G_2)=1$ if there is no prime number which divides both
  $\#G_1$ and $\#G_2$. 

\end{notation}

\begin{defns} Let $\call$ be a language, $T$ a complete theory. 
\begin{enumerate}
\item We say that the group $G$ is {\em geometrically represented in the theory $T$} if there exists $M_0\prec M \models T$ and $M_0 \subseteq A\subseteq B \subseteq M$,
$M$ such that   $\aut(B/A)\simeq G$, where $\aut(B/A)$ is the set of
permutations of $B$ which fix $A$ and preserve the truth value of all
$\call(A)$-formulas.  We say that a
{\em prime number $p$ is geometrically represented in $T$} if $p$ divides the
order of some finite group $G$ represented in $T$. 
%\item A field $F$ is {\em quasi-finite} if $G(F)\simeq \hat \zee$. \\
\item A field $F$ is {\em bounded} if for every integer $n$, $F$ has only
finitely many separable extensions of degree $n$. In this case we also
say that $G(F)$ is {\em bounded}. 
\item A field $F$ is {\em pseudo-algebraicallly closed}, henceforth
  abbreviated by {\em PAC}, if every absolutely irreducible variety
  defined over $F$ has an $F$-rational point. 
\item A field is {\em pseudo-finite} if it is PAC, perfect, and has
  exactly one extension of degree $n$ for each integer $n>1$. 

\end{enumerate}
\end{defns}
\begin{rem} (Folklore) Let $F$ be any field, $A$ a subfield of $F$, and assume that 
  $A=\dcl(A)$. Then $A^s\cap F$  is a Galois extension of
  $A$,  equals $\acl(A)$, and $\aut(\acl(A)/A)=\gal(A^s\cap F/A)$.  \\
Hence the finite groups $\aut(B/A)$ as above correspond to the
  finite quotients of $\gal(A^s\cap F/A)$. 

\end{rem}
\noindent 
Indeed, if $\alpha\in\acl(A)$, let $\alpha=\alpha_1,\alpha_2,\ldots, \alpha_n$
be the conjugates of $\alpha$ over $A$. Then the symmetric functions in
$\{\alpha_1,\ldots,\alpha_n\}$ are in $\dcl(A)=A$, i.e.: $\alpha$ satisfies a
monic separable polynomial with its coefficients in $A$ and $F$ contains all the
roots of this polynomial. This shows the first assertion and the second assertion is immediate. 

\para{\bf Properties of pseudo-finite fields and bounded PAC fields.} \\
We list some of the properties of these fields that we will use all the
time, often without reference. The language is the ordinary langage of
rings $\call=\{+,-,\cdot,0,1\}$, often expanded with parameters.
Pseudo-finite fields are the infinite models of the theory of finite
fields. They were studied by Ax in the 60's.  \\[0.1in]
An algebraic extension of a PAC field is PAC (Corollary 11.2.5 of \cite{FJ}). Theorem 20.3.3 of
\cite{FJ} (applied to $K=A$, $L=M=A^s\cap F$, $E=F=F$) gives the
following: \\[0.05in]
{\bf Fact 1}. {\em Let $F$ be a PAC field, $A$ a subfield of $F$ over which $F$ is
  separable, and assume that  $A$ has a Galois extension $C$ such
  that the restriction map $G(F)\to {\rm Gal}(C/A)$ is an
  isomorphism, and $C\cap F=A$. Let $B=A^s\cap F$; then $\aut_{\rm
    field}(B/A)=\aut(B/A)$. }\\ 
It suffices to notice that $CF=F^s$, and therefore also $CB=B^s$. So, if
$\varphi_0\in\aut_{\rm field}(B/A)$, extend $\varphi_0$ to
$\Phi_0\in\aut_{\rm field}(B^s/CA)$
by imposing $\Phi_0$ to be the identity on $C$. Then $\Phi_0$ induces
the identity on ${\rm Gal}(C/A)\simeq G(B)$. The result now follows
immediately from 20.3.3 in \cite{FJ}. It also has the following
consequence: \\[0.05in]
{\bf Fact 2}. {\em If $F_0\subset F$ are PAC fields of the same degree of
  imperfection, $F$ is separable over $F_0$, and the restriction
  map $G(F)\to G(F_0)$ is an isomorphism, then $F_0\prec F$.}\\[0.1in]
The following remark is totally folklore, but for want of a good
reference we will discuss it. \\[0.05in]
{\bf Fact 3}. {\em Let $F_0\prec F$ and assume that $G(F)$ is bounded. Then
the restriction
map $G(F)\to G(F_0)$ is an isomorphism.}\\
{}From $F_0\prec F$, it
follows immediately that $F$ is a regular extension of $F_0$, so that
the restriction map $G(F)\to G(F_0)$ is onto. Hence $G(F_0)$
 is bounded.  Fix an integer $n>1$,
and let $m(n)$ be the number of distinct separably algebraic extensions of $F_0$
of degree $n$. Then there is an $\call(F_0)$-sentence which expresses
this fact: that there are $m(n)$ distinct separably algebraic extensions
of $F_0$ of degree $n$, and that each separably algebraic extension
of degree  $n$ is contained in one of these. As $F_0\prec F$, $F$
satisfies the same 
sentence, and this implies that $F^s=F_0^s F$, and that the restriction
map $G(F)\to G(F_0)$ is an isomorphism.

\begin{lem}\vlabel{lem1} Let $F$ be a bounded field, and $A=\dcl(A)$ a
  subfield of $F$ containing an elementary substructure $F_0$ of $F$, 
  and let $B=A^s\cap F$. 
Then $G(A)\simeq G(F_0)\times
  \gal(B/A)$. 
\end{lem}

\prf  Because $G(F_0)$ is bounded and $F_0\prec F$, we know that
$F^s=F_0^sF$ and the fields $F_0^s$
and $F$ are linearly disjoint over $F_0$. Hence
$B^s=F_0^s B$, the fields $B$ and $AF_0^s$ are linearly disjoint over
$A$,  both are Galois extensions of $A$, and
therefore $G(A)=\gal(B^s/A)\simeq G(F_0)\times \gal(B/A)$.  

\begin{thm}\vlabel{koe} {\rm (Koenigsmann, Thm 3.3 in \cite{koe})}. Let $K$ be a field with
  $G(K)\simeq G_1\times G_2$. If a prime $p$ divides $(\# G_1,\# G_2)$, then there is
  a non-trivial Henselian valuation $v$ on $K$,  ${\rm char}(K)\neq p$,
  and ${\mathbb \mu}_{p^\infty}\subset K(\zeta_p)$. Furthermore, if $Kv$
  denotes the residue field of $v$ and $\pi:G(K)\to G(Kv)$ the canonical
  epimorphism, then $G(K)$ is torsion-free and $(\# \pi(G_1),\#
  \pi(G_2))=1$. 
\end{thm}

\begin{thm}\vlabel{thm1} Let $F$ be a field with bounded Galois group. Assume that $p$ is a
  prime number represented in ${\rm Th}(F)$ and  that $p$ divides $\#G(F)$. Then
  ${\rm char}(F)\neq p$, and $F(\zeta_p)$ contains ${\mathbb \mu}_{p^\infty}$. 
\end{thm}

\prf Let $F_0\prec F$, and $A$ a subfield of $F$ containing $F_0$, with
$A=\dcl(A)$. Let
$B=A^s\cap F$, and assume that $p$ divides $\#\gal(B/A)$, as well as
$\#G(F_0)$. By Lemma \ref{lem1}, we know that $G(A)\simeq G(F_0)\times
\gal(B/A)$. The result follows immediately from  Theorem \ref{koe}. 

\begin{thm}\vlabel{thm2} Let $F$ be a pseudo-finite field, [or more
  generally a bounded PAC field].  
Assume that for some  subfield  $A=\dcl(A)$ 
 of $F$ containing an elementary substructure of $F$, the group
 $G:=\aut(\acl(A)/A)$ is non-trivial. [Assume in addition
that all primes dividing $\#G$ divide $\#G(F)$]. \\
Then $G$ is abelian, and for any prime $p$ dividing $\#G$ we have $p\neq
{\rm char}(F)$ and 
$\mu_{p^\infty}\subset F$.
\end{thm}

\prf  Let $F_0\prec F$, and $A=\dcl(A)$ a subfield of $F$ containing
$F_0$, let
$B=A^s\cap F$, and assume that $p$ divides $\#\gal(B/A)$. By assumption,
$p$ divides 
$\#G(F_0)$, and by Lemma \ref{lem1}, $G(A)\simeq
\gal(B/A)\times G(F_0)$, with $p$ dividing the order of both
factors. Let $v$ be the Henselian valuation on $A$ given by Theorem \ref{koe},
and $\pi:G(A)\to G(Av)$ the corresponding epimorphism of Galois
groups. As $F_0$ is 
relatively algebraically closed in $A$, the valuation $v$ restricts to a
Henselian valuation on $F_0$; but because $F_0$ is PAC, the only Henselian
valuation on $F_0$ is the trivial valuation (\cite{FJ}, Cor 11.5.6). Hence $F_0\subseteq
Av$, and by Henselianity of $v$, $F_0^s\cap Av=F_0$. Hence the map $\pi$ is an isomorphism between $G(A)$ and
$G(F_0)$.  It follows that $\gal(BF_0^s/AF_0^s)$ is
contained in $\Ker(\pi)$, the inertia subgroup of $v$, and its order is
prime to the characteristic.  Hence $A^s$ is the composite of the purely
residual extension $AF_0^s$ of $A$,  and the totally ramified extension
$B$ of $A$. 
The characteristic of $F$ does not divide
$\# \gal(B/A)$, and this implies that $\gal(B/A)$ is abelian: indeed, by
Theorem 5.3.3 and \S~5.3 in 
\cite{EP}, we have 
$$\gal(B/A)\simeq \gal(BF_0^s/AF_0^s)\simeq 
{\rm Hom}(\Gamma(A^s)/\Gamma(AF_0^s)), {(Aw)^s}^\times),$$ where $w$
denotes the unique extension of $v$ to $A^s$, and $\Gamma(A^s)$,
$\Gamma(AF_0^s)$ the  value groups $w(A^s)$ and $w(AF_0^s)=v(A)$. \\
% The isomorphism is
% explicitly given by the following formula (see the beginning of \S5.3 in
% \cite{EP}): for $\delta\in \Gamma(B)$, and $a\in B$ such that
% $w(a)=\delta$, define $\psi(\si)(\delta)={\rm res}(\si(a)/a)$, where
% ${\rm res}$ is the residue map $B\to Bw$. Then $\psi(\si)$ does only
% depends on $\delta+\Gamma(A)$, and so defines a homomorphism from 
% $\Gamma(B)/\Gamma(A)$ into $(Bw)^s$. 
We also know that ${\mathbb \mu}_{p^\infty}\subset F(\zeta_p)$. Assume
first that $G(F)$ is abelian. Then so is 
$G(A)$, and therefore  any field between $A$ and $A^s$ is
a Galois  extension of $A$. In particular, because $p$ divide $\#G$, 
some element $\gamma\in v(A)$ is not divisible by $p$ in $v(A)$. Thus, if
$v(a)=\gamma$, then $a^{1/p}\in A^s$, and generates a Galois extension
of $A$: this implies that $\zeta_p\in A$, and by the above that
$\mu_{p^\infty}\subset F_0$. \\
Assume now that $G(F)$ is arbitrary, and that $\zeta_p\notin F_0$. Then
there is some $\si\in G(F)$ such that $\si(\zeta_p)\neq \zeta_p$, and
the subgroup  generated by $\si$ has order divisible by $p$ (here we use
that $p$ divides $\#G(F)$). Then the
restriction of $\si$ to $A^s$ commutes with all elements of
$\gal(A^s/F_0^sA)$, and so we may apply the previous result to the PAC
field $K$, subfield of $F^s$ fixed by $\si$, and its elementary
substructure $K_0$, subfield of $F_0^s$ fixed by $\si$, to deduce that
$\zeta_p\in K_0$, which contradicts our choice of $\si$.

\begin{cor} Let $F$ be a pseudo-finite field, or a bounded PAC field with $\# G(F)$ divisible by
  every prime number. Then every group
  represented in ${\rm Th}(F)$ is abelian. Furthermore, if $p$ is a
  prime represented in ${\rm Th}(F)$, then ${\mathbb
    \mu}_{p^\infty}\subset F$ and $p\neq {\rm char}(F)$. 
\end{cor}

\begin{cor} Let $F$ be a pseudo-finite field such that if $p$ is a prime
  number $\neq {\rm char}(F)$, then $\mu_{p^\infty}\not\subset F$. Then
definable closure and algebraic closure agree on subsets of $F$
containing an elementary substructure  of $F$. 
\end{cor}

\section{Other comments and remarks}
\para As was shown in Theorem 7  of \cite{BH}, if $F$ is a pseudofinite
field not of characteristic $p$ and containing $\mu_{p^\infty}$, then
every abelian $p$-group is represented in ${\rm Th}(F)$. Moreover, as the
class of groups represented in ${\rm Th}(F)$ is stable by direct product
(Remark 12 in \cite{BH}), it follows that which abelian groups are
represented in ${\rm Th}(F)$ is entirely determined by ${\rm char}(F)$
and by which $\mu_{p^\infty}$ are contained in $F$. \\
The proof given in \cite{BH} easily generalises to any perfect PAC 
field $F$, as they do have a notion of amalgamation over models, and
the construction did not use the pseudo-finiteness of $F$, only the fact
that it is PAC. We give here again the construction of a field with
absolute Galois group containing a cartesian product, it will be used in
the construction of example \ref{prop1}.

\para{\bf The construction}.  Let $F$ be a perfect field containing all
primitive roots of unity, and consider the field $K$ of generalized
power series $F^s((t^\rat))$ over $F^s$. Its members are formal sums
$\sum_{\gamma}a_\gamma t^\gamma$, with $\gamma\in\rat$, $a_\gamma\in F^s$, satisfying that $\{\gamma\mid
a_\gamma\neq 0\}$ is well-ordered. Then $K$ is algebraically closed. We
define an action of $G(F)$ on $K$ by setting $$\si(\sum_\gamma a_\gamma
t^\gamma) =\sum_\gamma \si(a_\gamma)t^\gamma.$$ So, the subfield of $K$
fixed by $G(F)$ coincides with $F((t^\rat))$. For each $n\in\nat$ not
divisible by the characteristic of $F$, choose a primitive $n$-th root
of unity $\zeta_n$, and choose them in a compatible way, i.e., such that
$\zeta_{nm}^m=\zeta_n$. Let $\si\in\aut(K)$ be defined by defining 
$\si(t^{1/n})=\zeta_nt^{1/n}$ for $n$ prime to the characteristic, and
if $q$ is a power of the characteristic, then $\si(t^{1/q})=t^{1/q}$; extend
$\si$ to the multiplicative group $t^{1/n}$, $n\in\zee$, and then to $K$
by setting
$$\si(\sum_\gamma a_\gamma
t^\gamma)=\sum _\gamma a_\gamma \si(t^\gamma).$$
Let $A$ be the subfield of $K$ fixed by $G(F)$ and by $\si$. Then
$G(A)\simeq G(F)\times \langle \si\rangle=G(F)\times \hat\zee$. 

\para{\bf Remark}. Let $F$ be a perfect PAC field, and let $A$ be the field
constructed above. So $A$ contains a copy of $F$ and is contained in
$F((t^\rat))$; as $F((t^\rat))$ is a regular extension of $F$, it
follows that $F$ has an elementary extension $F^*$ which contains
$B=A^s\cap F((t^\rat))$. Then $\aut(B/A)=\gal(B/A)\simeq \hat\zee$. This proof already
appeared in \cite{BH} (Thm 7).

\para {\bf Comment 1}. The proof of Lemma \ref{lem1} works exactly in
the same fashion as soon as the field $A$ contains enough information
about $G(F)$, more precisely: Assume $A$ contains $\acl(\emptyset)$, and
that for each finite extension $L$
of $F$, there is $\alpha$ such that $L=F(\alpha)$ and the minimal
polynomial of $\alpha$ over $F$ has its coefficients in $A=\dcl(A)$; then $A$
has a Galois extension $C$ which is linearly disjoint from $F$ over $A$,
and is such that $CF=F^s$. Then again one has $G(A)\simeq G(F)\times
\gal(A^s\cap F/A)$. The proof of Theorems \ref{thm1}  goes
through verbatim.\\
We were trying to weaken the hypotheses on $A$, and a natural weaker
assumption is to assume that $A$ contains a subfield $F_0$ such
that $\acl^{eq}(\emptyset)\subseteq \dcl^{eq}(F_0)$ and $\acl(F_0)=F_0$.
However the proof of Theorem \ref{thm2} used in an essential way the
fact that $F_0$ was PAC. The example below shows that this condition is
not sufficient. 

\para\vlabel{prop1} {\bf An example showing that the hypothesis of
  containing an elementary substructure is necessary}.  \\
Let $A_0$ be a field containing
  $\rat^{alg}$, and consider $A_0^{alg}((t^\rat))$; define  actions
  of $G(A_0)$ and of $\si$ on $A_0^{alg}((t^\rat))$ as
  above. Then $G(A_0((t)))\simeq G(A_0)\times\langle \si\rangle$. 
Let $F_0=\rat^{alg}((t))$,
the subfield of $\rat^{alg}((t^\rat))$ fixed by $\si$, and $A=A_0((t))$. Then
  $G(F_0)\simeq \hat\zee$,  and $A$ contains
  $F_0$. Furthermore, because $G(F_0)$
  is isomorphic to $\hat\zee$,  there is a pseudo-finite
  field $F$ which is a regular extension of $F_0$ (this follows easily from Thm 23.1.1 in \cite{FJ}), so that the
  restriction map $G(F)\to G(F_0)$ is an isomorphism. By Corollary~3.1 in \cite{H}, the theory of $F$
  eliminates imaginaries in the language augmented by constants for
  elements of $F_0$. As $F_0$ also contains
  $\acl(\emptyset)=\rat^{alg}$, it follows that
  $\acl^{eq}(\emptyset)\subset \dcl^{eq}(F_0)$. Furthermore, by standard
  results on pseudo-finite fields, $F$ has an elementary extension $F^*$
  which contains $A$ and is a regular extension of $B=A^{alg}\cap
  A_0^{alg}((t))$. Then $\gal(B/A)=\aut(B/A)\simeq G(A_0)$. \\
This shows that the hypothesis of $A$ containing an elementary
substructure of $F^*$ cannot be weakened to $A$ containing a  substructure
$F_0$ with $\acl^{eq}(\emptyset)\subset \dcl^{eq}(\emptyset)$ and $F_0=\acl(F_0)$.

\para {\bf Comment 2}. One can wonder what happens for a bounded PAC
field $F$ with $G(F)$ not divisible by all primes. If $S$ is the set of
prime numbers $\neq {\rm char}(F)$ and which do not divide $\#G(F)$, and
if $H$ is a projective $S$-group (i.e., the order of the finite
quotients of $S$ are products of members of $S$), then $G(F)\times H$
is a projective profinite group. Hence $F$ has a regular extension $K$
which is PAC and  
with $G(K)\simeq G(F)\times H$ (Thm 23.1.1 in \cite{FJ}). We may also impose, if the
characteristic is positive, that $K$ and $F$ have the same degree of
imperfection. As $K$ is a regular extension of $F$,
the restriction map $G(K)\to G(F)$ restricts to an isomorphism on
$G(F)\times (1)$, and sends $(1)\times H$ to $1$. Let $K_1$ be the
subfield of $K^s$ fixed by $G(F)\times (1)$. Then $K_1$ is PAC, and
because the restriction map $G(K_1)\to G(F)$ is an isomorphism, we have
$F\prec K_1$. If $A$ is the subfield of $K^s$ fixed by $G(F)\times H$,
then $A\subset F_1$, and $\gal(F_1/A)=\aut(F_1/A)\simeq H$. 

\para{\bf Comment 3}. Let $K$ be a  field,
$G=\aut(K(t)^{alg}/K(t))$, and $\si\in G$. Consider $G(\si)$ the
centralizer of $\si$ in $G$. Let $B$ be the subfield of $K(t)^{alg}$
fixed by $\si$, $F_0=K^{alg}\cap B$, and assume that $F_0$ is
pseudo-finite. Because $G(B)=\langle\si\rangle$ projects onto
$G(F_0)\simeq\hat\zee$, we have $G(B)\simeq \hat\zee$, and  $F_0$ has an elementary extension $F$ which is a
regular extension of $B$. We are interested in $\aut_{\rm
  field}(B/F_0(t))$; as $B\cap F_0^{alg}=F_0$, $B$ is linearly disjoint
from $F_0^{alg}(t)$ over $F_0(t)$, and therefore  $\aut_{\rm
  field}(B/F_0(t))=\aut(B/F_0(t))$, and its elements commute with
$\si$. \\
Let $H$ be a closed subgroup of
$G(\si)$ such that $H\cap \langle\si\rangle
=1$. Then Theorem \ref{thm2} tells us that $H$ is abelian, and that the
subfield $A$ of $B$ fixed by $H$ has a non-trivial Henselian valuation
$v$, which is trivial on $F_0$. Furthermore, if $p$ divides $\#H$, then $p\neq {\rm char}(F_0)$,
and $\mu_{p^\infty}\subset F_0$. We take the unique extension of $v$ to
$A^s$ (and also call it $v$); then the residue fields $Av$ and $Bv$ equal $F_0$,
and $(Av)^s =F_0^s$. Furthermore $H$ is procyclic, because
$\Gamma(A) \simeq \zee$, and $H\simeq {\rm Hom}(\rat/\zee,
{F_0^s}^\times)$. The  restriction of $v$ to $F_0(t)$ corresponds to a
point of $\proj^1(F_0)$ (because $Av=F_0$), i.e., either $v(t-a)=1$ for
some $a\in F_0$, or $v(t)=-1$. On the other hand, the field
$B$ can carry at most one Henselian valuation (see Thm
4.4.1 of \cite{EP}). It follows that $\aut_{\rm field}(B/F_0(t))$ is abelian,
procyclic. Hence  $G(\si)$ splits as
$\langle\si\rangle \times \langle \tau\rangle$. 
The result generalises
to any bounded PAC subfield $F_0$ of $K$, with exactly the same
reasoning. \\
This gives a partial
answer to Questions 15 and 16 of \cite{BH}. \\[0.05in]
Consider $K=\rat$, and endow $G(\rat(t))$ with the Haar measure. Then
the set $\{\tau\in G(\rat)\mid \rat^{alg}(\tau) \hbox{ is
  pseudofinite}\}$ has measure $1$, see Thm 18.6.1 in \cite{FJ}. Here
$\rat^{alg}(\tau)$ denotes the subfield of $(\rat^{alg})$ fixed
by $\si$. Moreover, it is easy to see that with probability $1$,
$\rat^{alg}(\si)$ does not contain $\mu_{p^{\infty}}$ for any prime
$p$. Hence, if $\si$ is any extension of $\tau$ to $\rat(t)^{alg}$, and
$B=\rat(t)^{alg}(\si)$, then $\aut(B/F_0(t))=1$.  \\[0.15in]
%
%
%\newpage
%\noindent
{\bf Addresses of the authors}\\
Bo\u{g}azi\c{c}i University\\
Faculty of Arts and Science \\
Department of Mathematics\\
34342, Bebek-Istanbul\\ 
Turkey\\[0.15in]
CNRS (UMR 8553) - Ecole Normale Sup\'erieure\\
45 rue d'Ulm\\
75230 Paris cedex 05\\
France

\end{document}